 \documentclass[12pt]{article}
 \usepackage[top=1in,bottom=1in,left=1.5in,right=1.5in]{geometry}
 \usepackage{amsmath,amssymb}
 \usepackage{latexsym}
 \usepackage[dvips]{pstricks} 
 \usepackage[dvips]{graphicx}

 \newcommand{\N}{\mathbb{N}}
 \newcommand{\R}{\mathbb{R}}
 \newcommand{\Z}{\mathbb{Z}}
 \renewcommand{\i}{\mathbf{i}}
 \renewcommand{\j}{\mathbf{j}}
 \newcommand{\e}{\mathbf{e}}
 
 \newcommand{\m}{\mathbf{m}}

 \newcommand{\T}{\mathbf{T}}

 \newcommand{\cA}{\mathcal{A}}
 \newcommand{\cB}{\mathcal{B}}

 \newcommand{\cG}{\mathcal{G}}

 \newcommand{\lan}{\langle}
 \newcommand{\ran}{\rangle}
 
 \newcommand{\hs}{\hspace*{\parindent}}
 \newcommand{\proof}{\hs \textbf{Proof.\ }}
 \newcommand{\tr}{\mathop{\mathrm{tr}}\nolimits}

 \newcommand{\qed}{\hspace*{\fill} $\Box$\\}
 \newcommand{\low}{\mathop{\mathrm{low}}\nolimits}
 \newcommand{\upp}{\mathop{\mathrm{upp}}\nolimits}

 \newcommand{\vol}{\mathrm{vol}}
 \newtheorem{theo}{\bfseries \hs Theorem}[section]
 \newtheorem{defn}[theo]{\bfseries \hs Definition}

 \newtheorem{corol}[theo]{\bfseries \hs Corollary}
 \newtheorem{con}[theo]{\bfseries \hs Conjecture}
 \newtheorem{example}[theo]{\bfseries \hs Example}
 \newtheorem{rem}[theo]{\bfseries \hs Remark}
 \numberwithin{equation}{section}

\begin{document}

 \title{On the validations of the \\ Asymptotic Matching Conjectures }

 \author{
 S. Friedland\thanks{Department of Mathematics, Statistics and Computer Science,
   University of Illinois at Chicago, Chicago, Illinois
   60607-7045, USA,
   and Berlin Mathematical School, Germany
  (friedlan@uic.edu).},
   E. Krop\thanks{Department of Mathematics, Statistics and Computer Science,
   University of Illinois at Chicago, Chicago, Illinois 60607-7045,
   USA (ekrop1@math.uic.edu).},
   P.H. Lundow\thanks{Department of Physics, AlbaNova University
 Center, KTH, SE-106 91 Stockholm, Sweden} and
 K. Markstr\"om\thanks{Deparment of
   Mathematics and Mathematical Statistics, Ume\aa University, SE-901
 87 Ume\aa, Sweden}}

   \date{April 5, 2008}

   \maketitle


 \begin{abstract}
      In this paper we review the asymptotic matching
      conjectures for $r$-regular bipartite graphs, and their
      connections in estimating the monomer-dimer entropies in
      $d$-dimensional integer lattice and Bethe lattices.
      We prove new rigorous upper and lower bounds for the
      monomer-dimer entropies, which support these conjectures.  We
      describe a general construction of infinite families of
      $r$-regular tori graphs and give algorithms for computing the
      monomer-dimer entropy of density $p$, for any $p\in [0,1]$, for
      these graphs.  Finally we use tori graphs to test the asymptotic
      matching conjectures for certain infinite $r$-regular bipartite
      graphs.  \\[\baselineskip] 2000 Mathematics Subject
      Classification: 05A15, 05A16, 05C70, 05C80,
      82B20\\[\baselineskip] Keywords and phrases: Matching and
      asymptotic growth of average matchings for $r$-regular bipartite
      graphs, monomer-dimer partitions and entropies.
 \end{abstract}


 \section{Introduction}

 The monomer-dimer covers of infinite graphs $G$, and in particular of
 the infinite graph induced by the lattice $\Z^d$, is one of the
 widely used models in statistical physics.  See for example
 \cite{Bax1, BoW, Ciu, Fis, FoR, FP, Gau, Ha1, Ha2, Ha3, HaM, HL, Kas,
 KeRaSi, Jer, L, Lun, Na2, Pau}.

 Let $G=(V,E)$ be an undirected graph with vertices $V$ and edges $E$.
 $G$ can be a finite or infinite graph.  A \emph{dimer} is a domino
 occupying an edge $e=(u,v)\in E$.  It can be viewed as two
 neighboring atoms occupying the vertices $u,v\in V$ and forging a
 bond between themselves.  A monomer is an atom occupying a vertex
 $w\in V$, which does not form a bond with any other vertex in $V$.  A
 \emph{monomer-dimer cover} of $G$ is a subset $E'$ of $E$ such that
 any two distinct edges $e,f\in E'$ do not have a common vertex.  Thus
 $E'$ describes all dimers in the corresponding monomer-dimer cover of
 $G$.  All vertices $V'\subset V$, which are not on any edge $e\in
 E'$, are the monomers of the monomer-dimer cover represented by $E'$.
 $E'$ is referred to here as a \emph{matching}.  $E'$ is called a
 \emph{perfect matching} if $V'=\emptyset$, i.e. all the vertices of
 $G$ are covered by the dimers.

 Consider first a finite graph $G=(V,E)$.  Then $E'$ is called an
 \emph{$l$-matching} if $\#E'=l$.  Note that $2l\le \#V$.  Let
 $\phi(l,G)\ge 0$ be the number of $l$-matchings in $G$ for any
 $l\in\Z_+$.  (Note that $\phi(0,G)=1$ and $\phi(l,G)=0$ if there are
 no $l$-matchings in $G$.  Assume also that $\phi(l,G)=0$ for a
 non-integer $l\ge 0$.)  Then the \emph{monomer-dimer} entropy of
 density $p$ of $G$ is defined as
 $$h_G(p)=\frac{\log \max(\phi(\lfloor\frac{p\cdot \#V}{2}\rfloor,G),1)}{\#V}
 \textrm{ for any } p\in[0,1].$$
 Let $\psi(x,G):=\sum_{l=0}^{\infty}
 \phi(l,G) x^l$ denote the matching generating polynomial of $G$.
 The \emph{pressure} of $G$ is defined as
 $$P_G(t):=\frac{\log\psi(e^{2t},G)}{\#V}.$$
 For an infinite graph $G$ the monomer-dimer entropy of density $p$ and the pressure
 $P_G(t)$ are defined by taking appropriate $\limsup$ on the
 finite sequences of graphs converging to $G$.  (See for details
 \S2.)

 We now consider the classical case in statistical physics: the
 lattice $\Z^d$, consisting of all $d$-dimensional vectors
 $\i=(i_1,\ldots,i_d)$ with integer coordinates.  (As usual we denote
 by $\Z,\Z_+,\N$ the set of integer, the set of nonnegative integers
 and the set of positive integers.)  Let
 $\e_k=(\delta_{k1},\ldots,\delta_{kd})$ be the unit vector in the
 direction of the coordinate $x_k$ for $k=1,\ldots,d$.  Then
 $G(\Z^d)=(V=\Z^d,E)$, where $(\i,\j)\in E\iff \j-\i=\pm\e_k$ for some
 $k\in [1,d]$.  Note that $G(\Z^d)$ is an infinite $2d$ regular graph.

 Let $h_d(p):=h_{G(\Z^d)}(p)$ for any $p\in [0,1]$ and
 $h_d:=\sup_{p\in[0,1]} h_d(p)$.  ($h_d$ and $\tilde h_d:=h_d(1)$ are
 called the $d$-monomer-dimer entropy and the $d$-dimer entropy
 respectively \cite{FP}.)  For $d=1$ it is known that \cite[\S4]{FP}:
 \begin{equation}\label{h1p}
      h_1(p)=\left( 1 - \frac{p}{2} \right) \log \left(1 - \frac{p}{2}
 \right) - \frac{p}{2} \log
      \frac{p}{2} - (1 - p)\log (1 - p), \quad p\in [0,1].
 \end{equation}

 The value of planar dimer entropy $h_2(1)$ was computed in \cite{Fis}
 and \cite{Kas} $$h_2(1)={1\over \pi}\sum_{q=0}^{\infty} {(-1)^q\over
 (2q+1)^2}= 0.29156090\ldots$$ The exact values of $h_2(p)$ for $p\in
 (0,1)$ and $h_d(p)$ for $d\ge 3, p\in (0,1]$ are unknown.  According
 to Jerrum \cite{Jer}, the computation of the matching generating
 polynomials
 of finite planar graphs in general is computationally intractable.
 (This fact does not rule out the possibility that $h_d(p)$ are
 computationally tractable for $d\ge 2$, however for $d\ge 3$ it seems
 that $h_3(1)$ and $h_3$ are hard to compute with high precision.)

 The properties of the entropy $h_d(p)$ for any $p\in [0,1]$ was
 studied by Hammersley and his collaborators in \cite{Ha1,Ha2, Ha3,
 HaM}.  It was shown in \cite{FP} that $h_d(p)$ can be obtained from
 the limits of certain tori graphs, which are bipartite and $2d$
 regular.  Using the proof of Tverberg's permanent conjecture, proved
 by the first name author \cite{Fr2}, the following lower bound was
 shown in \cite{FP}

 \begin{equation}\label{hfdineq}
 h_d(p)\ge fh_d(p):={\textstyle \frac{1}{2}} ( -p\log p - 2(1-p)\log(1-p) +
 p\log 2d -p )
 \end{equation}
  for any $p\in [0,1]$.

 Tverberg's permanent conjecture states that the minimum of the sum of
 all $l\times l$ permanental minors of $n\times n$ doubly stochastic
 matrices is achieved only at the flat matrix $J_n=(\frac{1}{n})$.  It
 is a generalization of the van der Waerden permanent conjecture for
 doubly stochastic matrices, which is the case $l=n$.  In \cite{Sch}
 Schrijver gave a lower bound on the number of perfect matchings for
 $r$-regular bipartite graphs.  It is an improvement of the lower
 bound implied by the van der Waerden permanent conjecture.
 Furthermore, this lower bound is asymptotically sharp.  Equivalently,
 one can think that Schrijver's lower bound gives asymptotically the
 number of perfect matchings in large random $r$-regular bipartite
 graph.

 In \cite{FKM} we stated a \emph{Lower Matching Conjecture}, referred
 here as LMC, for any $l$-matchings of $r$-regular bipartite graph.
 For $2l=\#V$ this conjecture is asymptotically equivalent to
 Schrijver lower bound for perfect matchings.  This lower bound can be
 viewed asymptotically as the number of $l$-matchings in a large
 random $r$-regular bipartite graph.  The LMC implies the \emph{Lower
 Asymptotic Matching Conjecture} stated in \S2, referred here as LAMC,
 yields the following conjecture.

 \begin{equation}\label{conjlbhdp}
      h_d(p)\ge gh_{2d}(p), \texttt{ for any } p\in (0,1] \texttt{ and
      } d\ge 2,
 \end{equation}
 where
 \begin{equation}\label{defghrp}
      gh_r(p):=\frac{1}{2}\big(p \log
      r -p\log p - 2(1-p)\log (1-p)
      +(r-p)\log (1 -\frac{p}{r})\big),
 \end{equation}
 for any integer $r \ge 2$.  Note that $h_1(p)=gh_2(p)$.  In a recent
 paper \cite{FG} the LAMC was proven for the
 sequence of densities $p=\frac{r}{r+s}, s=0,1,\ldots$, for any given $r\ge 2$.
 Hence (\ref{conjlbhdp}) holds for $p=\frac{2d}{2d+s}, s=0,1,\ldots$.  In
 particular $h_d(1)\ge gh_{2d}(1)$ for any $d\in\N$.  The inequality
 $h_3(1)\geq 0.440075$ is the best known lower bound.  A recent
 massive computation performed by the third named author in \cite{Lun}
 gives the best known upper bound $h_3(1) \leq 0.457547$.

 The conjectured lower bound (\ref{conjlbhdp}) yields a lower bound
 for the $d$-monomer-dimer entropy $h_d$.  In particular, the
 conjectured lower bound (\ref{conjlbhdp}) yields $h_3\geq
 0.784992989$.  The validity of (\ref{conjlbhdp}) for $d=3,
 p=\frac{6}{6+3}=\frac{2}{3}$ implies the best lower bound known
 $h_3\ge h_3(\frac{2}{3})\ge .7845241927$.  In this paper we give new
 lower bounds on $h_d(p)$ which yield the inequality $h_3\ge
 h_3(0.6814)\ge .7849602275$.  The numerical computations in \cite{FP}
 yield the best known upper bound $h_3\le 0.7862023450$.

 In \cite{FKM} we stated an \emph{Upper Matching Conjecture}, referred
 here as UMC. Namely, let $K_{r,r}$ be a complete bipartite graph on
 $2r$ vertices, where the degree of each vertex is $r$.  Denote by $q
 K_{r,r}$ be the graph consisting of $q$ copies of $K_{r,r}$.  Then
 the UMC claims that any $r$-regular bipartite graph $G$ on $2qr$
 vertices satisfies $\phi(l,G)\le \phi(l,q K_{r,r})$ for
 $l=0,\ldots,qr$.  We also have a corresponding \emph{Upper Asymptotic
 Matching Conjecture}, referred here as UAMC, which is slightly more
 technical to state.  (See \S6.)  For $r=2$ we proved these
 conjectures in \cite{FKM}.

 The main purpose of this paper is to give theoretical and numerical
 evidences on the LAMC and UAMC and their applications to the
 estimates of the monomer-dimer $p$-densities for $\Z^d$ and for the
 \emph{Bethe} lattices, i.e. $d$-regular infinite trees.  We believe
 that the computational and theoretical setting discussed in this
 paper are of interest by itself and to researchers in asymptotic
 combinatorics, which is widely used in statistical physics.

 We now outline briefly the main setting of our computations for the
 verification of the two asymptotic conjectures.  It is well known
 that the asymptotic growth of many configurations in statistical
 physics are given in terms of the spectral radius of the transfer
 matrix.  See for example \cite{FP}.  In this paper we construct
 infinite families $G_n=(V_n,E_n),n\in \N$ of $r$-regular bipartite
 graphs, which are coded by a specific incidence matrix $A\in
 \{0,1\}^{N\times N}$.  This sequence of graphs converges to an
 infinite $r$-regular graph $G$.  Using programs based on software
 developed by the third named author one obtains the transfer matrix
 $B(t) \in \R^{2^N\times 2^N}$, corresponding to the matching
 generating polynomial with the value $x=e^{2t}$.  Since the infinite
 tori graphs corresponds to a \emph{subshifts of finite type},
 abbreviated here as SOFT, one can compute the pressure function
 $P(t)$ in terms of the spectral radius $\rho(B(t))$.  This is well
 known to the experts, and we bring the proofs of these formulas in
 the paper for completeness, using the general techniques in
 \cite{FP2}.  (The properties of the \emph{pressure function}
 $P(t),t\in \R$, for multi-dimensional SOFT, as for example the
 monomer-dimer models in $\Z^d, d>1$, are studied in detail in
 \cite{FP2}.)  Then the monomer-dimer $p$-density $h_G(p)$ is computed
 by using $\rho(B(t))$ and its derivative.  (In this setting
 $p=p(t)$.)  We then compare $h_G(p(t))$ to the upper and lower bound
 given by the lower and upper asymptotic conjecture.

 We now briefly survey the contents of our paper.  In \S2 we discuss
 the monomer-dimer entropy $h_{\{G_n\}}(p)$ of density $p$ and the
 pressure function $P_{\{G_n\}}(t)$ for a sequence of finite graphs
 $\{G_n=(V_n,E_n)\}$ of bounded degrees such that $\#V_n\to\infty$.
 We define the function $\low_r(p), p\in [0,1]$ which gives the sharp
 inequality $h_{\{G_n\}}(p)\ge \low_r(p)$ for any sequence of
 bipartite $r$-regular graphs and any $p\in[0,1]$.  We state the LAMC,
 which is equivalent to the equality $\low_r=gh_r$.  Furthermore if
 the sequence $\{G_n\}$ is a sequence of random $r$-regular bipartite
 graphs we conjecture that $h_{\{G_n\}}(p)=gh_r(p)$ almost surely \cite{FKM}.
 In \S3 we use the recent verification of the LAMC for $p=\frac{r}{r+s},
 s=0,1,\ldots$ for any $r\ge 2$ to derive tight lower bounds on for
 $\low_r$.  In \S4 we discuss the applications of our results to Bethe
 lattices, i.e. infinite dimensional $r$-regular trees.  In \S5 we
 discuss the sequence of tori graphs, which are considered in
 \cite{FP} and \cite{FP2} to compute $h_2,h_3$ and $h_2(p)$.  We prove
 the thermodynamics formalisms for such graphs which gives the
 monomer-dimer entropy of density $p$ in terms of the pressure.  In
 \S6 we describe a fairly general construction of sequences of regular
 graphs, which includes the sequence of tori graphs.  In \S7 we
 describe the upper matching conjecture and its asymptotic version,
 called the upper asymptotic matching conjecture.  We give upper
 bounds for $h_{\{G_n\}}(p)$ for any sequence of bipartite $r$-regular
 graphs and show that in some regions these bounds are relatively
 close to the UAMC. In \S8 we describe our computational results,
 which support the conjectures stated in this paper.  In \S9 we
 identify an infinite graph with the maximal pressure among other
 infinite graphs in certain families of sequences described in \S6.

 \section{Entropies, Pressure and LAMC}
 We will now  define a limiting monomer-dimer density for a sequnce of
 bounded degree graphs.
 \begin{defn}\label{defconvig}
     Let $\{G_n\}, G_n=(V_n,E_n), n\in\N$ be a sequence of finite graphs,
     where multi edges are allowed, such that $\#V_n\to\infty$ and the
     degree of each vertex in $G_n$ is bounded by $d$ for $n\in N$.  For
     $p\in [0,1]$ we define $h_{\{G_n\}}(p)$, the monomer-dimer entropy of
     density $p$, as follows:
     \begin{eqnarray}\label{defpcon1}
     h_{\{G_n\}}(p)=\limsup_{n\to\infty}
     \frac{\log\phi(l_n,G_n)}{\#V_n},\\
     \texttt{ over all sequences }
     l_n\in \Z_+ \texttt{ satisfying }\lim_{n\to\infty} \frac{2l_n}{\#V_n}=p \in [0,1].
     \label{defpcon2}
     \end{eqnarray}
     $h_{\{G_n\}} (1)$ and $h_{\{G_n\}}:=\sup_{p\in [0,1]} h_{\{G_n\}}(p)$
     are called the dimer entropy of $\{G_n\}$,
     and the monomer-dimer entropy of $\{G_n\}$ respectively.
     For $t\in\R$ the pressure of $\{G_n\}$ is defined as
     \begin{equation}\label{defpres}
     P_{\{G_n\}}(t):=\limsup_{n\to\infty}
     \frac{\log\psi(e^{2t},G_n)}{\#V_n}.
     \end{equation}

     Let $G=(V,E)$ be an infinite graph, where multi edges are allowed.
     Assume that the maximal degree of vertices in $G$ is $d<\infty$.  A
     sequence of graphs $\{G_n\}, G_n=(V_n,E_n), \; V_n\subset V,
     n=1,2,\ldots, V_n$, where multi edges are allowed, converges to $G$
     if the following conditions hold:
     \begin{enumerate}
     \item
     $V_1\subset V_2 \subset\ldots$ are finite subsets of $V$ satisfying the
     condition $V=\cup_{n=1}^{\infty} V_n=V$.

     \item Each $G_n$ contains the induced subgraph $G(V_n):=(V_n,
     E(V_n))$ of $G$ on the set of vertices $V_n$, and the degree
     of each vertex in $G_n$ is at most $d$.

     \item Let $v\in V_n$ and assume that all neighbors of $v$ in
     $V$ are in $V_n$.  Then $E$ and $E_n$ have the same
     set of edges that contain $v$.
     \end{enumerate}

     Then
     $h_G(p):=h_{\{G_n\}}(p),\;h_G:=h_{\{G_n\}},\;P_G(t):=P_{\{G_n\}}(t)$.

 \end{defn}

 The above definition of entropy and pressure of an infinite graph $G$
 depends on the specific choice of the convergent sequence $\{G_n\}$
 to the infinite graph $G$.  For $G(\Z^d)$ one has a whole class of
 the sequences $\{G_n\}$, for which the resulting $h_{\{G_n\}}(p),
 P_{\{G_n\}}(t)$ is independent of the choice of the convergent
 sequence $\{G_n\}$ \cite{Ha1,Ha2, Ha3, HaM, FP2}.  In this case we
 denote by $h_d(p),P_d(t)$ the corresponding quantities.  For other
 infinite graphs $G$ discussed in this paper we choose a convenient
 convergent sequence $\{G_n\}$, and we do not discuss the
 corresponding class of sequences which yield the same entropy and
 pressure.

 The properties of the entropy $h_d(p)$ for any $p\in [0,1]$ was
 studied by Hammersley and his collaborators in \cite{Ha1,Ha2, Ha3,
 HaM}.  Let us mention two properties that are of interest in this
 context.  For any $m\in\N$ let $\lan
 m\ran:=\{1,\ldots,m\}=[1,m]\cap\Z$ be the set of integers between $1$
 and $m$.  For any $\m=(m_1,\ldots,m_d)\in\N^d$ let $\lan \m\ran:=\lan
 m_1\ran\times\ldots\times \lan m_d\ran\subset \N^d$ is the set of
 points in the lattice $\Z^d$ located in the box $[1,m_1]\times
 \ldots\times [1,m_d]$ in $\R^d$.  Denote by $\vol(\m):=\prod_{i=1}^d
 m_i$ the volume of the box $\lan \m\ran$.  Let $G(\m):=(\lan
 \m\ran,E(\m))$ be the subgraph of $G(\Z^d)$ induced by $\lan \m\ran$,
 i.e. $(\i,\j)\in E(\m) \iff \i,\j\in \lan \m\ran$ and $\i-\j=\pm
 \e_k$ for some $k\in \lan d\ran$.  Let
 $\m_n:=(m_{1,n},\ldots,m_{d,n}), n\in\N$ be a sequence of lattice
 points in $\N^d$, such that $\m_n\to\infty \iff m_{k,n}\to \infty$ as
 $n\to\infty$ for each $k=1,\ldots,d$.  Then for any sequence $l_n\in
 [0,\frac{\vol(\m_n)}{2}]\cap\N$ the following conditions hold:
 \begin{equation}\label{boxdefhdp}
      \lim_{n\to\infty} \frac{\log\phi(l_n,G(\m_n))}{\vol(\m_n)}=h_d(p)
      \texttt{ if } \m_n\to\infty \texttt{ and }
      \lim_{n\to\infty} \frac{2l_n}{\vol(\m_n)}=p\in[0,1] .
 \end{equation}
 The above characterization yields that $h_d(p)$ is a concave
 continuous function on $[0,1]$, see \cite{Ha1}.

 Let $T(\m):=(\lan\m\ran, \tilde E(\m))$ be the torus on $\lan
 \m\ran$.  Thus two vertices $\i,\j\in \lan\m\ran$ in $T(\m)$ are
 neighbors if $(\i,\j)\in E(\m)$, or for any $m_k>2$ the vertices
 $(i_1,\ldots,i_{k-1},1,i_{k+1},\ldots,i_d)$ and
 $(i_1,\ldots,i_{k-1},m_k,i_{k+1},\ldots,i_d)$ are adjacent for any
 $k\in \lan d\ran$ and \noindent
 $(i_1,\ldots,i_{k-1},i_{k+1},\ldots,i_d)\in \N^{d-1}$.  Clearly
 $\phi(l,T(\m))\ge \phi(l,G(\m))$ for any $l\in\Z_+$.  It was shown in
 \cite{FP} that the condition (\ref{boxdefhdp}) can be replaced by the
 corresponding condition on the torus:
 \begin{equation}\label{tordefhdp}
      \lim_{n\to\infty} \frac{\log\phi(l_n,T(\m_n))}{\vol(\m_n)}=h_d(p)
      \texttt{ if } \m_n\to\infty \texttt{ and }
      \lim_{n\to\infty} \frac{2l_n}{\vol(\m_n)}=p\in
      [0,1].
 \end{equation}
 (It is assumed that $l_n\in [0,\frac{\vol(\m_n)}{2}]\cap\N$.)
 More general, one can show that
 $h_d(p)=h_{G(\Z^d)}(p)$.

 There are several advantages of considering $T(\m)$ over $G(\m)$.
 Assume that $m_k >2$ for $k=1,\ldots,d$.  First, the graph $T(\m)$ is
 a $2d$-regular graph.  Second, the automorphism group of $T(\m)$ is
 quite big, which can be very well exploited, using the general
 method of \cite{Lun}.  See also\cite{LuMa}, and \cite{FP,FP2}
 for the computations of $h_d,\tilde h_d$  and
 $h_d(p)$ respectively.

 The fact that $T(2\m)$ is $2d$-regular bipartite graph was exploited
 in \cite{FP} to show (\ref{hfdineq}).  This lower bound is obtained
 by noting that if $G=(V,E)$ is an $r$-regular bipartite graph then
 $\phi(l,G)\ge f_r(l,\#V)$, where the function
 $f_r(l,2n):=\binom{n}{l}^2 l!\,\left(\frac{r}{n}\right)^l$ is
 determined from the proof of Tverberg's permanent conjecture
 \cite{Fr2}.

 The LMC stated in \cite{FKM} claims that $\phi(l,G)\ge g_r(l,\#V)\;(\ge
 f_r(l,\#V))$ for any $r$-regular bipartite graph, where.
 \begin{equation}\label{mmatchincon}
 g_r(l,2n)= {n\choose l}^2 (\frac{nr-l}{nr})^{rn-l} (\frac{lr}{n})^l.
 \end{equation}
 For $2l=\#V$ (\ref{mmatchincon}) is Schrijver's lower bound for
 perfect matchings in $r$-regular bipartite graphs on $2n$ vertices.
 The LAMC, which yields(\ref{conjlbhdp}), can
 be stated as follows:

 \begin{con}
     ( \textbf{The Lower Asymptotic Matching Conjecture}.)
     \label{lasmc} Let $\cG(2n,r)$ be the set of $r$-regular bipartite
     graphs on $2n$ vertices, possibly with multi edges.  For each
     $l\in [0,n]\cap\Z$ let $\mu(l,2n,r):=\min_{G\in
     \cG(2n,r)}\phi(l,G)$.  For $p\in [0,1]$ let $\low_r(p)$ be the
     infimum
     $\liminf_{n\to\infty} \frac{\log
     \mu(l_k,2n_k,r)}{2n_k}$
     over all sequences
     $0\le l_k\le n_k,k\in\N$ such that
    $lim_{k\to\infty} \frac{2l_k}{2n_k}=p$. Then
      \begin{equation}\label{lasmcp}
      \low_r(p)= gh_r(p).
      \end{equation}
 \end{con}

 The results in \cite{FKM} show that
 for a given $p\in [0,1]$ and
 $r\ge 2$, the above conjecture is equivalent to the statement that
 the number of $l$-matching in a random bipartite
 $r$-regular graph \emph{will} behave asymptotically as in
 Conjecture \ref{lasmc}.
 In particular, the random graphs minimize, in the asymptotical
 sense, the number of $l$-matchings in $r$-regular bipartite
 graphs.

 It is shown in \cite{FKM} that the LMC and LAMC hold for $r=2$.
 Furthermore, the cycle $C_{2n}$ on $2n$ vertices satisfies the
 inequality $\phi(l,C_{2n})\le \phi(l,G)$ for any $G\in\cG(2n,2)$.
 Hence
 $$\lim_{k\to\infty} \frac{\log\phi(l_k,C_{2n_k})}{2n_k}=gh_2(p)=h_1(p)$$
 for any sequence $0\le l_k\le n_k, k\in\N$ satisfying the assumptions of
 Conjecture \ref{lasmc}.

 In a recent paper \cite[Theorem 5.6]{FG} the following results were
 proven:
 \begin{theo}\label{FGthm}
     The Lower Asymptotic Matching Conjecture holds for the following
     corresponding sequence of densities $p=\frac{r}{r+s},
     s=0,1,\ldots$, for any given $r$.  In particular (\ref{conjlbhdp})
     holds for $p=\frac{2d}{2d+s}, s=0,1,\ldots$.
 \end{theo}
 This can be extended to give a bound for all $p$ in the following way.
 \begin{defn}\label{defntilg}
     For $2\le r\in N$ let $ghl_r(p), p\in [0,1]$ be the following
     function.
     \begin{itemize}
     \item
     $ghl_r(\frac{r}{r+s})=gh_r(\frac{r}{r+s})$ for
     $s=0,1,\ldots$.

     \item $ghl_r(p)$ is linear on the interval
     $[\frac{r}{r+s+1},\frac{r}{r+s}]$ for $s=0,1,\ldots$.

     \item $ghl_r(0)=0$.

     \end{itemize}

 \end{defn}

 The concavity of $h_d(p)$ and Theorem \ref{FGthm} yields:
 \begin{equation}\label{lowbndhd}
     h_d(p)\ge ghl_{2d}(p), \textrm{ for all } p\in [0,1] \textrm{ and }
     d=2,3,\ldots.
 \end{equation}
 In the next section we improve substantially these lower bounds.

 In \cite[Figure 1]{FG} are plotted the graph of $gh_{4}(p)$, the
 graph corresponding to UAMC and the 19 values of the $h_2(p)$
 computed by Baxter \cite{Bax1}.  (Baxter's computations are based on
 sophisticated heuristical arguments.  His computations were recently
 verified by rigorous mathematical methods in \cite{FP2}.)  It turns
 out that Baxter's values are very close to the values of $gh_{4}(p)$.

\begin{figure}[tbp]
      \includegraphics*[width=14cm,height=8cm]{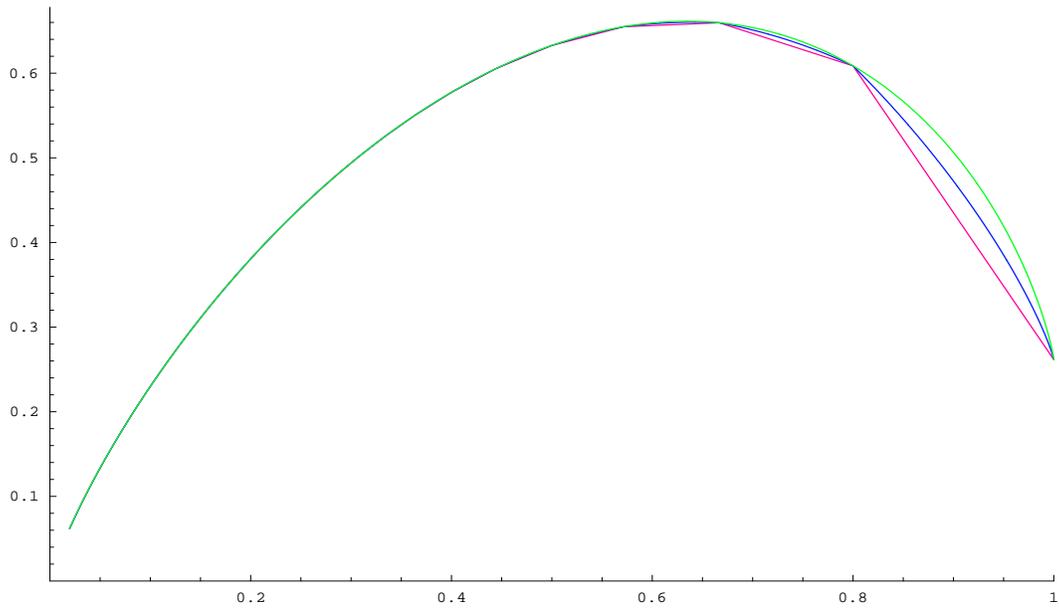}
      \caption{$ghl_4$-red, $\low_{4,1}$-blue, $gh_4$-green, color online}
      \label{fig:f2}
 \end{figure}

 In Figure \ref{fig:f2} we show the graphs of
 $ghl_r(p)$, $\low_{r,1}(p)$, a lower bound for $\low_r(p)$ given in
 the next section, and $gh_r(p)$ for $r=4$.  Note that the
 differences of the three graphs are relatively large on the first
 interval from the right $[\frac{r}{r+1},1]$, slightly less on the
 second interval from the right $[\frac{r}{r+2},\frac{r}{r+1}]$, and
 ignorable from the fourth interval to the right
 $[\frac{r}{r+3},\frac{r}{r+2}]$.  We notice that the differences
 between the functions $ghl_r(p),\low_{r,1}(p),gh_r(p)$ decrease as
 $r$ increases.  (This observation applies also for the values $r=3,6$
 which are not plotted here.)

 \section{Lower bounds for $\low_r(p)$}
 In this section we give a lower bound for the function $\low_r(p)$,
 which is defined in Conjecture \ref{lasmc}.
 \begin{theo}\label{wconv}
     The function $\low_r(p)+\frac{1}{2}(p\log p +(1-p)\log(1-p))$ is
     concave.
 \end{theo}
 \proof  Let $G\in \cG(2n,r)$.
 Consider the polynomial $\theta(x):=(-x)^n\psi(-\frac{1}{x},G)$.  Since
 $\phi(l,G)>0$ for $l=0,\ldots,n$ it follows that $\theta(x)$ has
 $n$ complex nonzero roots.
 It is well known \cite{HL} that $\theta(x)$ has only positive
 roots.  The Newton inequalities, see e.g \cite{Nicu}, yield
 \begin{equation}\label{newtin}
     \frac{\phi(l-1,G)}{{n\choose l-1}}\frac{\phi(l+1,G)}{{n\choose
     l+1}}\le \big(\frac{\phi(l,G)}{{n\choose l}}\big)^2, l=1,\ldots,n-1.
 \end{equation}
 Let $G_{l,2n,r}\in \cG(2n,r)$ be an $r$-regular graph for which the
 equality $\mu(l,2n,r)=\phi(l,G_{l,2n,r})$ holds.  (\ref{newtin}) and
 the minimal characterizations of $\mu(k,2n,r), r=0, \ldots,n$ yields
 \begin{equation}\label{newtin1}
     \frac{\mu(l-1,2n,r)}{{n\choose l-1}}\frac{\mu(l+1,2n,r)}{{n\choose
     l+1}}\le \big(\frac{\mu(l,2n,r))}{{n\choose l}}\big)^2, l=1,\ldots,n-1.
 \end{equation}
 This is equivalent to the statement that the sequence $$a_{l,2n,r}:=
 \log \frac{\phi(l,G)}{{n\choose l}}\ge 0,\quad  l=0,\ldots,n$$
 is a concave sequence.   Let $\alpha(x,2n,r)$ be a piecewise linear
 function on $[0,1]$ defined as follows:
 \begin{itemize}
     \item
     $\alpha(\frac{l}{n},2n,r)=\frac{a_{l,2n,r}}{2n}, \quad
     l=0,\ldots, n$.

     \item $\alpha(\frac{l}{n},2n,r)$ is linear function on the
     interval $[\frac{l}{n}, \frac{l+1}{n}]$ for $l=0,\ldots, n-1$.
 \end{itemize}
 The concavity of the sequence $a_{l,2n,r}, l=0,\ldots,n$ is
 equivalent to the concavity of $\alpha(x,2n,r)$.  Let
 $\{n_k\}_{k=1}^{\infty}$ be an increasing sequence of positive
 integers.  Let $l_k\in [0,n_k]\cap\Z, k\in\N$ be a sequence
 satisfying $\lim_{k\to\infty} \frac{l_k}{n_k}=p\in [0,1]$.  Use
 Stirling's formula to deduce that
 \begin{equation}\label{constirf}
     \lim_{k\to\infty} \frac{\log{n_k \choose
     l_k}}{2n_k}=-\frac{1}{2}\big(p\log p+(1-p)\log(1-p)\big).
 \end{equation}
 It is straightforward to show that
 \begin{equation}\label{limchlow}
     \liminf_{n\to\infty} \alpha(x,2n,r)=\low_r(x)+
     \frac{1}{2}\big(p\log p+(1-p)\log(1-p)\big).
 \end{equation}
 Since each $\alpha(\cdot,2n,r)$ is concave it follows that
 $\low_r(p)+
 \frac{1}{2}\big(p\log p+(1-p)\log(1-p)\big)$ is concave.

 \qed

 The arguments of the proof of the above Theorem combined with the
 definition of $h_d(p)=h_{G(\Z^d)}(p)$, implies a stronger concavity
 result than given in \cite{Ha1}.

 \begin{corol}\label{sconchd}
     Let $h_d(p), p\in [0,1]$ be the
     monomer-dimer entropy of density $p$ for the graph $G(\Z^d)$.
     Then $h_d(p)+
     \frac{1}{2}\big(p\log p+(1-p)\log(1-p)\big)$ is concave.
 \end{corol}

 Theorems \ref{FGthm} and \ref{wconv} yields:

 \begin{corol}\label{lblow}
     Let $\low_{r,1}(p), p\in
     [0,1]$ be defined as follows.
     \begin{itemize}
     \item
     $\low_{r,1}(\frac{r}{r+s})=gh_r(\frac{r}{r+s})$ for
     $s=0,1,\ldots$.

     \item $\low_{r,1}(p)+\frac{1}{2}(p\log p +(1-p)\log(1-p))$ is linear on the interval
     $[\frac{r}{r+s+1},\frac{r}{r+s}]$ for $s=0,1,\ldots$.

     \item $\low_{r,1}(0)=0$.

     \end{itemize}

     Then $\low_r(p)\ge \low_{r,1}(p)$ for any $p\in [0,1]$.

 \end{corol}

 Figure \ref{fig:f2} shows the position of the graphs
 $ghl_r(p)\le \low_{r,1}(p) \le gh_r(p)$ for $r=4$.

 We now give a different lower bound for $\low_r(p)$ using
 \cite[Theorem 5.6]{FG}.

 \begin{theo}\label{lblow2}
     Let $\low_{r,2}(p), p\in [0,1]$ be defined as follows.
     \begin{itemize}
     \item
     $\low_{r,2}(\frac{r}{r+s})=gh_r(\frac{r}{r+s})$ for
     $s=0,1,\ldots$.
     \item For $p \in (\frac{r}{r+1},1)$
     $\low_{r,2}(p)$ the maximum between the two following numbers
     \begin{eqnarray*}
         &&\frac{p}{2} \big(\log r +(r-1)\log(1-\frac{1}{r})\big)-\frac{1}{2}
         \big(p\log p + (1-p)\log (1-p)\big)\\
         &&\textrm{and}\\
         &&\frac{p}{2} \log r-p\log p - (1-p)\log (1-p) +
         \frac{r}{2}
         \log(1-\frac{1}{r+1}).
     \end{eqnarray*}

     \item For $p \in (\frac{r}{r+s+1},\frac{r}{r+s})$
     $\low_{r,2}(p)$ the maximum between the two following numbers
     \begin{eqnarray*}
         &&\frac{p}{2} \log r+\frac{1}{2}\left(-p\log p - 2(1-p)\log (1-p)\right) +\\
         &&\frac{1}{2}\left(
         (r+s-1)\log(1-\frac{1}{r+s})-(s-1+p)\log(1-\frac{1-p}{s})
         \right)\\
         &&\textrm{ and }\\
         &&\frac{p}{2} \log r+\frac{1}{2}\left(-p\log p - 2(1-p)\log (1-p)\right) +\\
         &&\frac{1}{2}\left(
         (r+s)\log(1-\frac{1}{r+s+1})-(s+p)\log(1-\frac{1-p}{s+1})
         \right).
     \end{eqnarray*}

     for $s=1,\ldots$.
     \item $\low_{r,2}(0)=0$.
     \end{itemize}

     Then $\low_r(p)\ge\low_{r,2}(p)$ for any $p\in [0,1]$.

 \end{theo}
 \proof   \cite[Theorem 5.6]{FG} states
 \begin{eqnarray}\label{lasrppbp}
 \liminf_{n \to \infty} \frac{\log \phi(l_k,G_k)}{2n_k} \geq
 \frac{p}{2}\log r+
 \frac{1}{2}\left(-p\log p - 2(1-p)\log (1-p)\right) +\\
 \frac{1}{2}\left(
 (r+j-1)\log(1-\frac{1}{r+j})-(j-1+p)\log(1-\frac{1-p}{j})
 \right),\nonumber
 \end{eqnarray}
 for any sequence of $G_k\in \cG(2n_k,r)$ where $\lim_{k\to\infty}
 \frac{l_k}{n_k}=p$ and any $j\in \N$.
 This implies the the inequality $\low_r(p)\ge\low_{r,2}(p)$ where
 $p\in [\frac{r}{r+s+1},\frac{r}{r+s}]$ and $s\in\N$.
 (Choose $j=s,s+1$.)

 In the interval $[\frac{r}{r+1},1]$
 the second inequality follows from the above
 inequality for $j=1$.
 For the first inequality
 we use the arguments of Theorem
 \ref{wconv}. Combine the arithmetic-geometric inequality and Schrijver's
 inequality \cite{Sch} to deduce
 $$\big(\frac{\phi(l,G)}{{n \choose l}}\big)^{\frac{1}{l}}\ge
 \phi(n,G)^{\frac{1}{n}}\ge \frac{(r-1)^{r-1}}{r^{r-2}}.$$

 Taking logarithm of both sides, dividing by $2n$ and using any
 sequence satisfying (\ref{defpcon2}) we deduce the inequality
 \begin{equation}\label{frgur1}
 \low_r(p)\ge \frac{p}{2} \big(\log r +(r-1)\log(1-\frac{1}{r})\big)-\frac{1}{2}
 \big(p\log p + (1-p)\log (1-p)\big)
 \end{equation}
 for all $p\in [0,1]$.
 \qed

\begin{figure}[tbp]
      \includegraphics*[width=14cm,height=8cm]{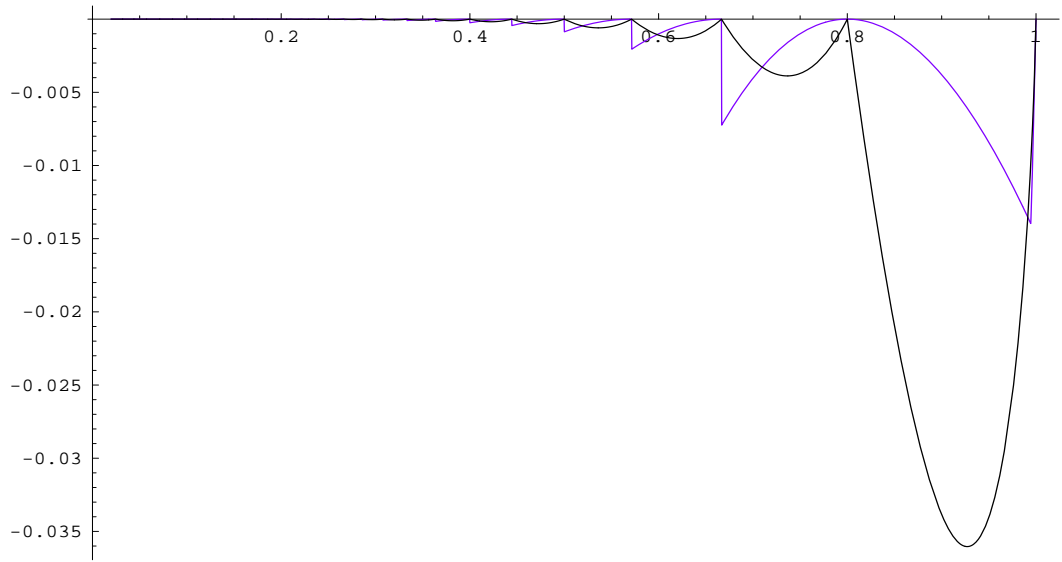}
      \caption{$\low_{4,1}-gh_4$-black, $\low_{4,2}-gh_4$-blue, color online}
      \label{fig:f8}
 \end{figure}

 It turns out that for many of the values of $p$, the lower bound
 $\low_{r,2}(p)$ a better lower bound than $\low_{r,1}(p)$, and it is very close
 to the function $gh_r(p)$.
 Figure \ref{fig:f8} compares the differences
 $\low_{r,1}(p)-gh_r(p)$, plotted in black, and
 $\low_{r,2}(p)-gh_r(p)$, plotted in blue, for $r=4$.

 From this graph and the graphs for $r=3,6$, which are not
 plotted here, we conclude that the errors
 $gh_r(p)-ghl_r(p), gh_r(p)-\low_{r,1}(p), gh_r(p)-\low_{r,2}(p)$
 are decreasing monotonically with $r$.

 Since $\low_{6,2}(0.6814)=.7849602275$ we deduce the inequality
 $$h_3\ge h_3(0.6814)\ge \low_{6}(0.6814)\ge
 \low_{6,2}(0.6814)=.7849602275$$
 given in \S1.  Combine Corollary \ref{lblow} with Theorem \ref{lblow2}
 to deduce

 \begin{corol}\label{lblowm}  For $3\le r\in\N$ and $p\in [0,1]$
 $$\low_r(p)\ge \max(\low_{r,1}(p),\low_{r,2}(p)).$$
 \end{corol}

 \section{Monomer-dimer densities for Bethe lattices}

 Let $\T(r)$ be an infinite $r$-regular tree.  Recall that $\T(3)$ is
 known the \emph{Bethe lattice}.  Clearly, each $\T(r)$ is bipartite.
 For each $\T(r), r\ge 2$ we construct a convergent sequence $\{G_n(r)\}$
 in sense of Definition \ref{defconvig}.

 Fix a vertex $O$ in $\T(r)$ and consider all vertices in $\T(r)$
 whose distance from $O$ is $n\ge 1$.  Then the number of such
 vertices is $r (r-1)^{n-1}$.  Let $A_1, \ldots,A_r$ be the $r$
 vertices of distance $1$ from $0$.  The number of vertices in $\T(r)$
 whose distance from $O$ is exactly $n$ is divided to $r$ classes
 $\cA_{i,n}, i=1,\ldots,r$, where the points in $\cA_{i,n}$ have
 distance $n-1$ from $A_i$.  Let $X_n=\cup_{i=1}^{r-1}\cA_{i,n}$.
 Note that $\#X_n=\#\cA_{r,n+1}=(r-1)^n$.  Let
 $V_n:=\{O\}\cup_{i,j=1}^{r,n}\cA_{i,j}\cup A_{r,n+1}$ Let
 $H_n(r)=(X_n\cup \cA_{r,n+1}, F_n)$ be an \emph{arbitrary} $r-1$ regular
 bipartite graph with the two classes of vertices $X_n,\cA_{r,n+1}$.
 Let $G_n(r)=(V_n,E_n)$ where $E_n$ are the union of the edge set in
 the induced graph $\T(r)(V_n)$ and the set $F$.  Note that $G_n(r)$
 is $r$-regular and bipartite.  Then $G_n(r), n=1,\ldots,$ converges
 to $\T(r)$.

 Note that $\T(2)$ is isomorphic to the integer lattice $\Z$, and
 $G_n(2)$ is a cycle of length $2(n+1)$ for $n\in\N$.  Hence
 \begin{equation}\label{hdt2}
     h_{G(\Z)}(p)=h_{\{G_n(2)\}}(p)=gh_2(p) \textrm{ for all } p\in [0,1].
 \end{equation}
 In \cite{FKM} it is shown that $\low_2(p)=gh_2(p)$.

 Using the definition $\low_r(p)$ and the fact that $G_n(r)$ are
 $r$-regular and bipartite we obtain.
 \begin{corol}\label{hdtr}
     Let $2<r\in\N$.  Then
     \begin{equation}\label{hdtrin}
     h_{\{G_n(r)\}}(p)  \ge \low_r(p) \textrm{ for all } p\in [0,1].
     \end{equation}
 \end{corol}
 For more complex models, like the Ising model, it is known that this
 kind of limit is sensitive to the exact limiting sequence of graphs
 \cite{Hagg}.
 It is an interesting problem if equality holds in the above
 inequality for some choices of random graphs $H_n(r)\in
 \cG(2(r-1)^{n},r-1), n\in\N$.

 \section{An example of sequence of tori}

 We first discuss
 a sequence of graphs that give the lower and upper bounds for $h_d$
 and $h_d(1)$ for the graph $G(\Z^d)$ considered in \cite{FP}.
 Assume that the dimension $d>1$.  Let
 $\m':=(m_1,\ldots,m_{d-1})\in \N^{d-1}$ be fixed and assume that
 $m_i>3$ for $i=1,\ldots,d-1$. Consider the sequence of
 $d$-dimensional tori $T((\m',n))=(V_n,E_n), n=3,4,\ldots$. Each
 torus is a $2d$-regular graph.  If $m_1,\ldots,m_{d-1}$ and $n$
 are even then $T((\m',n))$ is bipartite. The vertex set of
 $T((\m',n))$ is the set $V_n:=\lan (\m',n)\ran$. $V_n$ can be
 viewed as composed of $n$ layers of vertices $\lan \m'\ran$.  The
 edges between all vertices $\lan \m'\ran$ in each level $k$ are
 given as in the $d-1$ dimensional torus $T(\m')$.  The other edges
 of $T((\m',n))$ are going from level $k$ to level $k+1$ for
 $k=1,\ldots,n$, where the level $n+1$ is identified with the level
 $1$.  (We also identify level $0$ with the level $n$.)  The rule
 for the edges between the level $k$ and the level $k+1$ is
 independent of $k$.  Thus the vertices $(\i',k)$ and $(\j',k+1)$
 in $V_n$ are adjacent if and only if $\i'=\j'$.  The adjacency
 matrix between the two vertices $\i'$ in the level $k$ and $\j'$
 in the level $k+1$ is given by the $0-1$ matrix
 $A(\m'):=(a_{\i'\j'})_{\i',\j'\in \lan\m'\ran}$, which is an
 identity matrix of order $N=\vol(\m')$.  For any square matrix
 $A\in \R^{n\times n}$ we denote by $\tr A$ and $\rho(A)$ the trace
 and the spectral radius of $A$ respectively.

 Let us recall first the computation of the monomer-dimer entropy
 $h_d$ given in \cite{FP}.  The entries of the transfer matrix
 $B(\m')=(b_{ST})_{S,T\subset \lan \m'\ran}$ are indexed by two
 subsets of $S,T$ of $\lan \m'\ran$.  (These subsets may be empty.)
 First $b_{ST}=0$ if $S\cap T\ne \emptyset$.  Second assume that
 $S\cap T=\emptyset$ then $b_{ST}$ counts the number of the
 monomer-dimer covers of the subgraph of $T(\m')$ induced by the
 set vertices $\lan\m'\ran\backslash (S\cup T)$.  Note that any
 subgraph of $T(\m')$ induced by a set $U\subset \lan \m'\ran$ can
 be covered by monomers.  Hence $b_{ST}\ge 1$.  (If $S\cup
 T=\lan\m\ran$ then $b_{ST}=1$.)  It is not hard to see that the
 product of $n$ terms $b_{S_1S_2}b_{S_2S_3}\ldots
 b_{S_{n-1}S_{n}}b_{S_{n}S_1}$ corresponds to all monomer-dimer
 covers of $G_n=T((\m',n))$ with the following conditions.  For
 each level $k=1,\ldots,n$ the dimers going from the level $k$ to
 $k-1$ are located at the set $S_k$ and the dimers going from the
 level $k$ to the level $k+1$ are located at the set $S_{k+1}$. Let
 $\Phi(G_n):=\sum_{l=0}^{\infty} \phi(l,G_n)$ number of all
 possible monomer-dimer covering $T((\m',n))$.  Then
 $\tr B(\m')^{n}=\Phi(G_n)$.
 It is shown in \cite{FP}
 \begin{eqnarray}\label{limmdgr}
      &&\lim_{n\to\infty}
      \frac{\log\Phi(G_n)}{\#V_n}=\frac{\log\rho(B(\m'))}{\vol(\m')},\;
      h_d=\lim_{\m'\to\infty} \frac{\log\rho(B(\m'))}{\vol(\m')},\\
      &&\texttt{ and } h_d\le \frac{\log\rho(B(2\m'))}{\vol(2\m')}
      \texttt{ for any } \m'\in\N^{d-1}.\nonumber
 \end{eqnarray}
 (Here $2\m':=(2m_1,\ldots,2m_{d-1})$.)
 The lower bounds for $h_d$ are also expressed in terms
 of linear combinations of certain $\log\rho(B(\m'))$ corresponding
 to different values of $\m'$.

 Let $T((\m',\Z))$ be an infinite graph given by the set of
 vertices $\lan \m'\ran \times \Z$ and the following set of edges
 $\tilde E(\m',\Z)$.  $((\i',p),(\j',q))\in \tilde E(\m',\Z)$ if
 either $p=q$ and $(\i',\j')\in \tilde E(\m')$ or $|p-q|=1$ and
 $\i'=\j'$.  Thus the sequence of graphs $G_n:=T((\m',n)),
 n=3,4,\ldots$ converges to $G:=T((\m',\Z))$.  Let $h_G(p):=h_{\{G_n\}}(p)$ be
 defined by (\ref{defpcon1}-\ref{defpcon2}). We now show how to
 compute $h_G(p)$ using the pressure function.

 Let $S,T$ be two disjoint subsets of $\lan \m'\ran$. Let
 $E'\subset \tilde E(\m')$ be an $l$-matching of $T(\m')$ so that
 each edge $(u,v)\in E'$ represents a dimer occupying two adjacent
 vertices in $T(\m')$ located in $\lan \m'\ran\backslash (S\cup
 T)$.  To this matching we correspond a monomial $x^l$.  Let
 $c_{ST}(x)$ be the sum of all such monomials. $c_{ST}(x)$ is the
 matching polynomial for the graph $T(\m', S,T)$, the
 subgraph of $T(\m')$ induced by the subset of vertices $\lan
 \m'\ran\backslash (S\cup T)$.  We let $c_{ST}(x)=0$ if $S,T\subset
 \lan \m'\ran$ and $S\cap T\ne \emptyset$.  Then
 $b_{ST}=c_{ST}(1)$. Let $b_{ST}(t):=c_{ST}(e^{2t})e^{(\#S+\#T)t}$
 and $B(\m,t):=(b_{ST}(t))_{S,T\subset \lan \m'\ran}$.
 The arguments in \cite{FP} that show that $\Phi(G_n)=\psi(1,G_n)$
 yield the equality $\tr B(\m,t)^n=\psi(e^{2t},G_n)$.
 The definition (\ref{defpres}) of
 pressure $P(t)$ and the arguments in \cite{FP} for the equality
 (\ref{limmdgr}) imply
 \begin{equation}\label{defpresm}
 P(t):=P_{\{G_n\}}(t)=\frac{\log\rho(B(\m,t))}{\vol(\m')}, \quad t \in \R.
 \end{equation}
 (We suppressed the dependence of $P(t)$ on $\m'$.)

 The following results are known, e.g. \cite{Bax1, FP2},
 and we bring their proof for completeness.

 \begin{theo}\label{presprop}
 Let $P(t):=P_{\{G_n\}}$ be defined by (\ref{defpresm}).
 Then $P$ is a smooth increasing convex function on $\R$.
 Furthermore
 \begin{equation}\label{rangept}
 \lim_{t\to -\infty} P'(t) =0, \quad \lim_{t\to\infty} P'(t) =1.
 \end{equation}
 Let $p(t):=P'(t), t\in \R$.  Then $p(t)$ is an increasing function
 such that $p([-\infty,\infty])=[0,1]$.  Furthermore
 \begin{equation}\label{hgpeq}
      h_{\{G_n\}}(P'(t))=P(t)-tP'(t)  \texttt{ for any }t\in\R,
 \end{equation}
 and $h_{\{G_n\}}(p)$ is a continuous concave function on $[0,1]$,
 which is smooth in $(0,1)$.
 \end{theo}

 \proof
 The well known result \cite{Kin} yields $P(t)$ is a convex
 function of $t$.  Since $\m'$ is fixed we let $B(t):=B(\m',t)$
 and $N=\vol(\m')$.
 As $B(t)$ is irreducible, and the nonzero entries of $B(t)$
 are increasing function on $\R$, it follows that
 $P(t)$ increases.  Since $\rho(B(t))$ is a simple root positive of the
 characteristic polynomial it follows that $\rho(B(t))$ and
 $P(t)$ is an analytic function in some open domain containing
 $\R$.
 \begin{equation}\label{bmtform}
 B(t)=\sum_{i=0}^{N} e^{it}B_i,\textrm{ each } B_i \textrm{ is
 nonnegative}.
 \end{equation}
 $B_0=(b_{ST,0})_{S,T\subset \lan \m'\ran}$ corresponds to monomer configurations
 of . I.e.
 $B_0$ is a diagonal matrix with one nonzero entry $b_{\emptyset
 \emptyset,0}=1$.  The matrix $B_{N}$ corresponds to
 the tiling of $T((\m',n))$ by dimers.  Hence
 $P(t)=\log\rho(B(t))$ is not a linear function.
 The analyticity of $P(t)$ yields that $P''(t)$ may have only finite
 number of zeros on any closed interval $[a,b]$.  The convexity of
 $P(t)$ implies that $P''\ge 0$ on $\R$.  Hence $P''$ is positive on
 any $[a,b]$ except a finite number of points.  Thus $p(t)=P'(t)$
 increases on $\R$ and $P(t)$ is strictly convex on $\R$.
 Let $x=e^{t}$.  Then the $\tilde B(x):=B(\log x)$ is a polynomial in $x$.
 Since $\rho(\tilde B(0))=\rho(B_0)=1$ is a simple root of the
 characteristic polynomial of $B_0$ it follows that
 $\rho(\tilde B(x))$ is analytic in some disk $|x|<\epsilon$,
 such that $|\rho(\tilde B(x))|>0$ in this disk.
 Hence the branch $\log \rho(\tilde B(x)), \log \rho(\tilde
 B(0))=0$ is analytic in this disk and has Taylor expansion.
 The same statement holds for the derivative of  $\log \rho(\tilde B(x))$.
 Substitute $x=e^t$ to deduce that
 $P(t)=\frac{\log \rho(\tilde B(e^t))}{N}$ and its derivative have
 convergent series in $x=e^{t}$ for $t<-R$, for some $R>>1$.
 This implies the first equality in (\ref{rangept}).

 Observe that
 $$P(t)=t+ \frac{\log\rho(\hat  B(t))}{N},\quad \hat  B(t):=\sum_{i=0}^{N}
 e^{(i-N)t}B_i.$$
 Hence $P'(t)=1+ (\frac{\log\rho(\hat  B(t))}{N})'$.
 The arguments above for the first equality in (\ref{rangept}) imply that the
 second equality in (\ref{rangept}).

 We now show the inequality
 \begin{equation}\label{plowbnd}
 P(t)\ge pt+h_G(p) \textrm{ for all } p\in [0,1].
 \end{equation}
 Let $l_n\in [0,n]\cap \Z, n\in\N$ be  sequence satisfying (\ref{defpcon2}).
 Then
 \begin{eqnarray*}
 \tr B(t)^n =\psi(e^{2t},G_n)\ge
 e^{2l_nt}\phi(l_n,G_n)\Rightarrow\\
 \frac{\log \tr B(t)^n}{N n}\ge \frac{2l_nt}{N n}
 +\frac{\log\phi(l_n,G_n)}{N n}.
 \end{eqnarray*}
 Recall that $\#V_n$, the number of vertices in $G_n$ is
 $N n$.  Use the definition of $P(t)$ and
 (\ref{defpcon2}) to deduce
 $$P(t)\ge pt
 +\limsup_{n\to\infty}\frac{\log\phi(l_n,G_n)}{\#V_n}.$$
 Use the definition (\ref{defpcon1}-\ref{defpcon2}) of $h_G(p):=h_{\{G_n\}}(p)$
 to deduce the inequality (\ref{plowbnd}).
 It is straightforward to show that $h_{G}(p)$ upper semicontinuous
 on $[0,1]$.

 We now show that for each $t\in\R$ there exists $p(t)\in [0,1]$ such that:
 \begin{equation}\label{pupbnd}
 P(t)\le p(t)t+h_G(p(t)).
 \end{equation}
 Let $l_n(t)\in [0,\frac{n N}{2}]\cap \Z]$ satisfy
 $$e^{2lt}\phi(l,G_n)\le e^{2l_n(t)t}\phi(l_n(t),G_n), \textrm{ for } l=0,\ldots
 \lfloor \frac{n N}{2}\rfloor.$$
 Hence
 \begin{equation}\label{pupbnd1}
 \phi(e^{2t},G_n)\le
 \frac{n N}{2}e^{2l_n(t)t}\phi(l_n(t),G_n).
 \end{equation}

 Take a subsequence $n_k, k\in\N$ such that
 $\lim_{k\to\infty} \frac{\log \psi(e^{2t},G_{n_k})}{n_k
 N}=P(t)$.  Choose a subsequence $j_k, k\in\N$ of $n_k,
 k\in\N$, such that $\lim_{k\to\infty}
 \frac{l_{j_k}(t)}{j_k N}=p(t)\in [0,1]$.
 Take the logarithm of the inequality (\ref{pupbnd1}) and divide by
 $n N$.  Let $n=j_k$ and let $k\to\infty$.
 The definition (\ref{defpcon1}-\ref{defpcon2}) of $h_G(p)$
 yields the inequality (\ref{pupbnd}).

 The inequalities (\ref{plowbnd}) and (\ref{pupbnd}) yield the
 equality $P(t)=p(t)t+h_G(p(t))$.  Moreover
 (\ref{plowbnd}) yields $P(y)$ lies above the line $p(t)y+h_G(p(t)$,
 which intersect $P(y)$ at the point $y=t$.  Hence $p(t)=P'(t)$ and
 $h_G(p(t))=P(t)-p(t)t$.  I.e. (\ref{hgpeq}) holds.  Since $P'$
 increasing and analytic the implicit function theorem yields that
 $t=Q(p)$ is analytic in $p\in (0,1)$.  Hence $h_G(p)$ is analytic
 on $(0,1)$.  Observe that $-h_G(p)$ is the Legendre function
 corresponding to a smooth strictly convex function $P(t)$
 \cite{Roc}.  Hence $h_G(p)$ is concave on $(0,1)$.  Our arguments
 yield that $h_G(p)$ is continuous on $[0,1]$.  Hence $h_G(p)$ is a concave
 function on $[0,1]$.  \qed

 \begin{rem}\label{termdynfor} Let
 $h_G(p):=h_{\{G_n\}}(p)$ and $P(t):=P_{\{G_n\}}(t)$ be given as in Definition
 \ref{defconvig}.
 Theorem \ref{presprop} applies to $h_G(p)$ and $P(t)$ in the
 following cases:
 \begin{enumerate}
 \item  There exists a nonnegative irreducible matrix $B(t)$
 of the form (\ref{bmtform}) such that
 \begin{itemize}
 \item $\#V_n = n N, \;n\in\N$.
 \item $\psi(e^{2t},G_n)=\tr B(t)^{n},\;n\in\N$.
 \item
 $\rho(B_0)=1$ and $\rho(B_N)$
 are positive simple roots the characteristic polynomials of
 $B_0$ and $B_N$ respectively.
 \end{itemize}

 2.  $G_n$ is a disjoint union of $n$ copies of a finite graph
 $H=(W,F)$ which has a perfect matching.  Then
 $P(t)=\frac{\log\psi(e^{2t},H)}{\#W}$.
 \end{enumerate}
 \end{rem}

 \section{A construction of sequences of graphs}\label{seq}

 We now generalize the construction in the previous section to a
 general construction of a sequence of regular graphs.  Let
 $F=(U,D)$ be an undirected graph with the set of vertices $U$ and
 the set of edges $D$.  For $n\ge 2$ let $G_n:=(V_n,E_n)$ be the
 following graph. $V_n=U\times \lan n\ran$, i.e. we can view $V_n$
 consisting of $n$ copies of $U$ arranged in the $n$ layers
 $(U,1),(U,2),\ldots,(U,n)$. We let $(U,0):=(U,n), (U,n+1):=(U,1)$.
 Then
 \begin{enumerate}
      \item For any $u,v\in U$  and $k\in \lan n \ran$
      $((u,k),(v,k))\in E_n \iff (u,v)\in D$.

      \item Any other edges of $E_n$ are between the vertices $(U,k)$
      and $(U,k+1)$ for $k=1,\ldots,n$.

      \item Let $A=(a_{uv})_{u,v\in U}$ be a given nonzero $0-1$ matrix.
      Then for each $k\in \lan n\ran$ $((u,k),(v,k+1))\in E_n \iff
      a_{uv}=1$.  We call $A$ the \emph{connection} matrix.

      \item
      For any two subsets $S,T\subset U$, ($S,T$ may be empty), let
      $\tilde a_{ST}\in\Z_+$ be defined as follows.  If $\# S\ne \#T$
      then $\tilde a_{ST}=0$.  Assume that $\# S=\#T$.  Let $\cB(S,T)$
      be the set of all bijections $\beta:S\to T$.  Then $\tilde
      a_{\emptyset \emptyset}=1$ and $\tilde
      a_{ST}=\sum_{\beta\in\cB(S,T)} \prod_{s\in S} a_{s\beta(s)}$ for
      $\#S=\#T\ge 1$.

      Thus $\tilde a_{ST}$ is the number of perfect matchings in the
      subgraph of the bipartite graph on the set of vertices $(U,1)\cup
      (U,2)$, and the set edges $E\subset (U,1)\times (U,2)$ given by
      $A$, and induced by the subset of vertices $(S,1)\cup (T,2)$.  Let
      $\tilde A:=(\tilde a_{ST})_{S,T\subset U}$ be a $2^{\#U}\times
      2^{\#U}$ matrix with nonnegative integer entries.

      \item
      For any two disjoint subsets $S,T\subset U$, let $c_{ST}(x)$ be
      the matching generating polynomial of the subgraph of $F$ induced by
      the set of vertices $U\backslash (S\cup T)$.  For non-disjoint
      subset $S,T\subset U$ let $c_{ST}(x)=0$.  Let
      $M(t):=(c_{ST}(e^{2t})e^{(\#S+\#T)t})_{S,T\subset U}$ and
      $B(t):=M(t)\tilde A$ be
      $2^{\#U}\times 2^{\#U}$ nonnegative matrices for any $t\in \R$.
      Then $\log\rho(B(t))$ is a continuous convex function on
      $\R$.  If $B(1)$ is an irreducible matrix then
      $\log\rho(B(t))$ is an analytic function on $\R$.
      (See arguments of the proof of Theorem \ref{defpresm}.)
 \end{enumerate}

 Then the sequence $G_n, n=2,\ldots$ has the following properties:
 \begin{itemize}
      \item
      If $F$ is connected then each $G_n$ is connected.

      \item
      Assume that $F$ is bipartite, where $U=U_1\cup U_2, D\subset
      U_1\times U_2$.  Suppose that the edges between the two
      consecutive levels of vertices $(U,k)$ and $(U,k+1)$ are either
      between $(U_i,k)$ and $(U_i,k+1)$ for $i=1,2$ or between $(U_i,k)$
      and $(U_{i+1},k+1)$ for $i=1,2$.  ($U_3:=U_1$.)  If $n$ is even
      then $G_n$ is bipartite.

      \item
      Assume that $F$ is $p$-regular.  Assume that the matrix $A$ has
      $q$ $1$'s in each row and column.  Then $G_n$ is $p+2q$-regular
      graph.

      \item
      Assume that $F$ is $p$-regular bipartite.  Let $U=U_1\cup U_2,
      D\subset U_1\times U_2$ and $n$ is even.  Assume that the matrix
      $A$ has the following properties.  Each row indexed by $u\in U_1$
      and each column indexed by $v\in U_2$ has $q$ $1$'s, and each row
      indexed by $v\in U_2$ and each column indexed by $u\in U_1$ has
      $q-1$ $1$'s.  ($q\in\N$.)  Then $G_n$ is $p+2q-1$ regular.

      \item
      Then sequence of graphs $G_n, n=2,3,\ldots$ converges to the
      infinite graph $G=(V,E)$, where $V=F\times \Z$.  The edges $E$ are
      either between the two vertices on the same level $(U,k), k\in
      \Z$, determined by $D$, or between the vertices of two consecutive
      levels $(U,k)$ and $(U,k+1)$, given by the incidence matrix $A$ in
      the way described above.

      \item
      $P(t):=\frac{\log\rho(B(t))}{\#U}$ is the pressure of $G$.
      Assume that $B(1)$ is an irreducible matrix.  Let $h_G(p)$
      be defined by (\ref{defpcon1}-\ref{defpcon2}).  Then (\ref{hgpeq})
      holds.  (See Remark \ref{termdynfor}.)

 \end{itemize}

 In the example of $G_n=T((\m',n)), n=3,4,\ldots$, discussed in the
 previous section, we have that $U=T(\m')$ and $A$ is the identity
 matrix $I$.  Hence $\tilde A$ is also the identity matrix.

 \section{The upper matching conjecture}

 For $r\ge 2$ let $K_{r,r}$ be a complete bipartite graph on $2r$
 vertices, where each vertex has degree $r$.  Then
 \begin{equation}\label{matkrr}
      \phi(l,K_{r,r})={r\choose l}^2 l!, \;l\in\Z_+, \texttt{ and }
      \psi(x,K_{r,r})=\sum_{l=0}^r {r\choose l}^2 l! x^l.
 \end{equation}

 \begin{con}\label{upmatcon} ( \textbf{The upper matching conjecture}.)
      Let $G=(V,E)$ be a finite bipartite
      regular $r$-regular graph on $2qr$ vertices where $2\le
      q,r\in\N$.  Let $q K_{r,r}$ be the graph consisting of $q$ copies
      of $K_{r,r}$.
      Then  $\phi(l,G)\le
      \phi(l,q K_{r,r})$ for $l=0,\ldots,qr$.
 \end{con}

 In \cite{FKM} we proved the above conjecture for $r=2$.  We also
 showed that for $r=2$ $\phi(l,G)\le \phi(l,q K_{2,2})$ for any $2$
 regular graph $G$ on $4q$ vertices.  ($G$ does not have to be
 bipartite.)  It is plausible that in the above conjecture one can drop
 the assumption that $G$ is bipartite.
 For $l=0,1$ the above conjecture is trivial.  For $l=qr$ the above
 conjecture follows from the Minc conjecture proved by Bregman
 \cite{Bre}.

 Let $K(r)$ be an infinite countable
 union of $K_{r,r}$.  Let $h_{K(r)}(p)$ be defined as in
 (\ref{defpcon1}-\ref{defpcon2}) where $G_n=n K_{r,r}, n\in\N$.
 Let $G_n=(V_n,E_n), n\in \N$ be a sequence of $r$ regular bipartite
 graphs, where $\#V_n \to\infty$.  Let $h_{\{G_n\}}(p)$ be defined as in
 (\ref{defpcon1}-\ref{defpcon2}).  Assume
 for simplicity of the exposition that $\#V_n=2q_n r$.  Then Conjecture
 \ref{upmatcon} yields $\phi(l,G_n)\le \phi(l,{q_n} K_{r,r})$ for $n\in
 \N$.  Hence the UMC yields the AUMC:
 $h_G(p)\le h_{K(r)}(p)$ for any $p\in [0,1]$.

 We use the pressure $P_{K(r)}(t)$, as pointed in
 Remark \ref{termdynfor}, to compute $h_{K(r)}(p)$.
 Clearly the
 matching generating polynomial of $q K_{r,r}$ is $\psi(x,K_{r,r})^q$.  Hence
 \begin{equation}\label{presK}
      P_{K(r)}(t)=\frac{\log \sum_{l=0}^r {r\choose l}^2 l! e^{2lt}}{2r},
      \quad t\in\R.
 \end{equation}
 This formula follows also from the results of the previous section,
 where $F=K_{r,r}$ and the incidence matrix $A$ between two levels
 $(U,1)$ and $(U,2)$ is the zero matrix.  Then $\rho(B(t)\tilde
 A)=\psi(e^{2t},K_{r,r})$. (\ref{hgpeq}) yields
 \begin{eqnarray}\label{hKp}
      &&h_{K(r)}(p(t))=\frac{\log \sum_{l=0}^r {r\choose l}^2 l!
      e^{2lt}}{2r}-\frac{t\sum_{l=0}^r {r\choose l}^2 l! (2l)e^{2lt}}
      {2r \sum_{l=0}^r {r\choose l}^2 l! e^{2lt}},\\
      && \texttt{where }p(t)=\frac{\sum_{l=0}^r {r\choose l}^2 l! (2l)e^{2lt}}
      { 2r\sum_{l=0}^r {r\choose l}^2 l! e^{2lt}} \texttt{ and } t\in\R.
      \nonumber
 \end{eqnarray}

 \begin{con}\label{upasmatcon}
      ( \textbf{The upper asymptotic matching conjecture}.)
      Let $G_k=(V_k,E_k), k\in \N$ be a sequence of $r$ regular bipartite
      graphs, where $\#V_k \to\infty$.   Let $h_{\{G_k\}}(p)$ be defined as
      in (\ref{defpcon1}-\ref{defpcon2}).  Let $h_{K(r)}(p)$ be defined by
      (\ref{hKp}).  Then $h_{\{G_k\}}(p)\le h_{K(r)}(p)$ for any $p\in [0,1]$.
 \end{con}

 It is plausible to assume that Conjecture \ref{upasmatcon} holds
 under the assumption that each $G_n$ is an $r$-regular graph.

 \begin{theo}\label{upmatchpbd}  Let $2\le r\in\N$ and assume that
 $G_n=(V_n,E_n), n\in N$ is a sequence of $r$-regular bipartite
 graphs such that $\#V_n\to\infty$.   Let
 $h_{\{G_n\}}(p), p\in [0,1]$ be defined by
 (\ref{defpcon1}-\ref{defpcon2}).
 Then
 \begin{equation}\label{upmatchpbd1}
 h_{\{G_n\}}(p)\le \min(\upp_{r,1}(p),\upp_{r,2}(p)) \textrm{ for all
 } p\in [0,1],
 \end{equation}
 where
 \begin{eqnarray}
 &&\upp_{r,1}(p):=\frac{p\log r!}{2r} -\frac{1}{2}p\log p -(1-p)\log(1-p)
 \label{defupp1},\\
 &&\upp_{r,2}(p):=\frac{p\log r}{2} -\frac{1}{2}\big(p\log p +
 (1-p)\log(1-p)\big) \label{defupp2}.
 \end{eqnarray}

 \end{theo}

 \proof  We claim that for any $G=(V_1\cup V_2,E)\in \cG(2m,r)$ the
 following inequality holds
 \begin{equation}\label{fupbnd1}
 \phi(l,G)\le {m\choose
 l}\frac{(m!)^{\frac{m-l}{m}}}{(m-l)!} (r!)^{\frac{l}{r}}.
 \end{equation}
 Indeed, let $U_1\subset V_1$ be a subset of
 cardinality $l$.  Consider the induced bipartite graph
 $G(U_1\cup V_2)$.  Then $G(U_1\cup V_2)$ is also an induced
 subgraph of the graph $G'=(V_1\cup V_2,E')$, where
 the induced subgraph $G'(V_1\backslash U_1
 \cup V_2)$ is the complete bipartite graph on $V_1\backslash
 U_1\cup V_2$.  It is straightforward to see that
 $\phi(l,G')=(m-l)!\phi(l,G(U_1\cup V_2)$.
 The Bregman inequality \cite{Bre} yields
 $\phi(l,G')\le (m!)^{\frac{m-l}{m}} (r!)^{\frac{l}{r}}$.
 Since the number of choices of $U_1$ is ${m\choose l}$ we deduce
 the inequality (\ref{fupbnd1}).
 Let $G_n=(V_n,E_n)\in \cG(2m_n,r), n\in\N$.
 Let $l_n\in [0,m_n]\cap \Z$ be a sequence satisfying (\ref{defpcon2}).  Take the logarithm of
 the (\ref{fupbnd1}) divide by $\#V_n=2m_n$ and let $n\to\infty$ to
 obtain that $h_{\{G_n\}}(d)\le \upp_{r,1}(p)$.

 Our next inequality is
 \begin{equation}\label{fupbnd2}
 \phi(l,G)\le {m \choose l} r^{l}.
 \end{equation}
 Let $U_1\subset V_1, \#U_1=l$.  Since each vertex of $G$ has degree
 $r$ it follows that each vertex in $U_1$ that $\phi(l,G(U_1\cup
 V_1))\le r^l$.  As one have $n \choose l$ choices of the set $U_1$
 we obtain (\ref{fupbnd2}).
 The above arguments imply that
 $h_{\{G_n\}}(p)\le \upp_{r,2}(p)$.  Hence (\ref{upmatchpbd1})
 holds.  \qed

\begin{figure}[tbp]
      \includegraphics*[width=14cm,height=8cm]{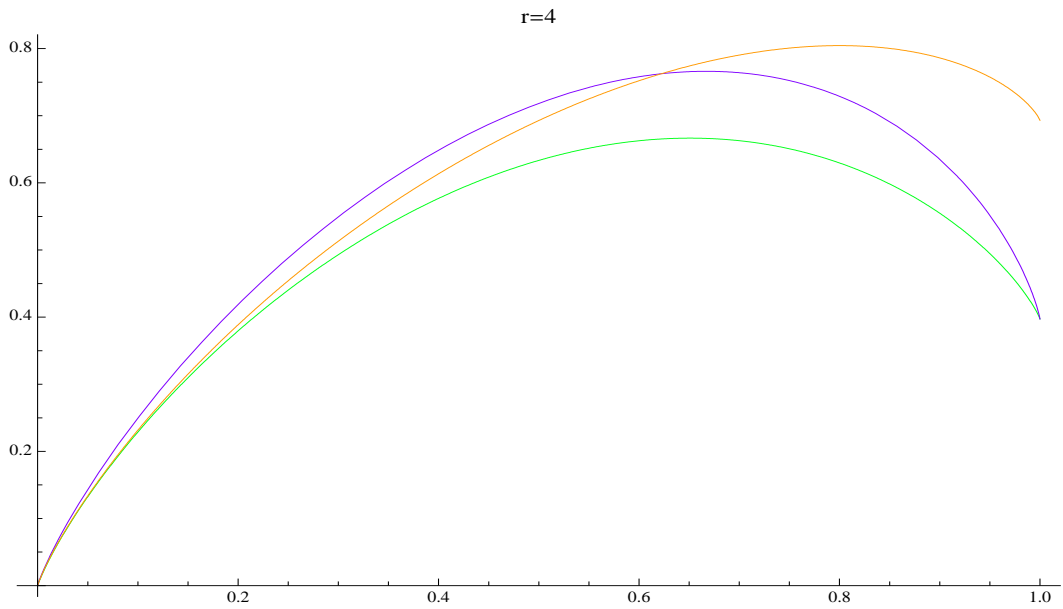}
      \caption{$h_{K(4)}$-green, $\upp_{4,1}$-blue, $\upp_{4,2}$-orange, color online}
      \label{fig:f11}
 \end{figure}

 Figure \ref{fig:f11} gives the plot of
 $h_{K(r)},\upp_{r,1},\upp_{r,2}$ for $r=4$.
 From this graph and the corresponding graphs
 for $r=3,8$, which is not plotted here, we see that
 $\min(\upp_{r,1}(p),\upp_{r,1}(p))-h_{K(r)}(p)$ decreases with $r$.
 Moreover the intersection point of the graphs $\upp_{r,1}$ and
 $\upp_{r,1}$ moves to the left as $r$ increases.

 \section{Computational results}

 We have checked the asymptotic matching conjectures for several
 families like those described in Section \ref{seq}. In each case
 we choose $U$ to be a cycle $C_l$ of length $l$ for several values
 of $l$ and then varied the connection matrix $A$. In each case we
 used the described transfer matrix method to compute the entropy
 for  several values of $p$ and then compared with the conjectured
 bounds. In all cases the conjectures were found to hold.

 In order to test the conjectured lower bound for a given choice of
 $U$ and $A$ we first constructed the transfer matrix $B(t)$ for
 the given graph. Given $B(t)$ we can directly compute $P(t)$ from
 the maximum eigenvalue as in (\ref{defpresm}). Next we computed
 $P'(t)$, using the equality
 $$\rho'(t)=\eta_{1}^T\left(\frac{d}{dt}B(t)\right) \eta_{2},$$
 where $\eta_{1}^T$ and $\eta_{2}$ are the left and right
 eigenvectors of $B(t)$, normalized by the condition
 $\eta_{1}^T\eta_{2}=1$.  (This is a standard variational formula,
 e.g. \cite{Fr1}.)
 From these values we now compute
 $h_G(p(t))$ using (\ref{hgpeq}). So for each value of $t$ we get a
 pair $(p(t),h_G(p(t)))$ telling us the asymptotic pressure
 $h_G(p(t))$ at the density $p(t)=P'(t)$. To make all computations
 exact we chose $e^t$ to be rational numbers, which yielded
 rational values for all matrix entries.

\begin{figure}[tbp]
     \includegraphics*[width=14cm,height=8cm]{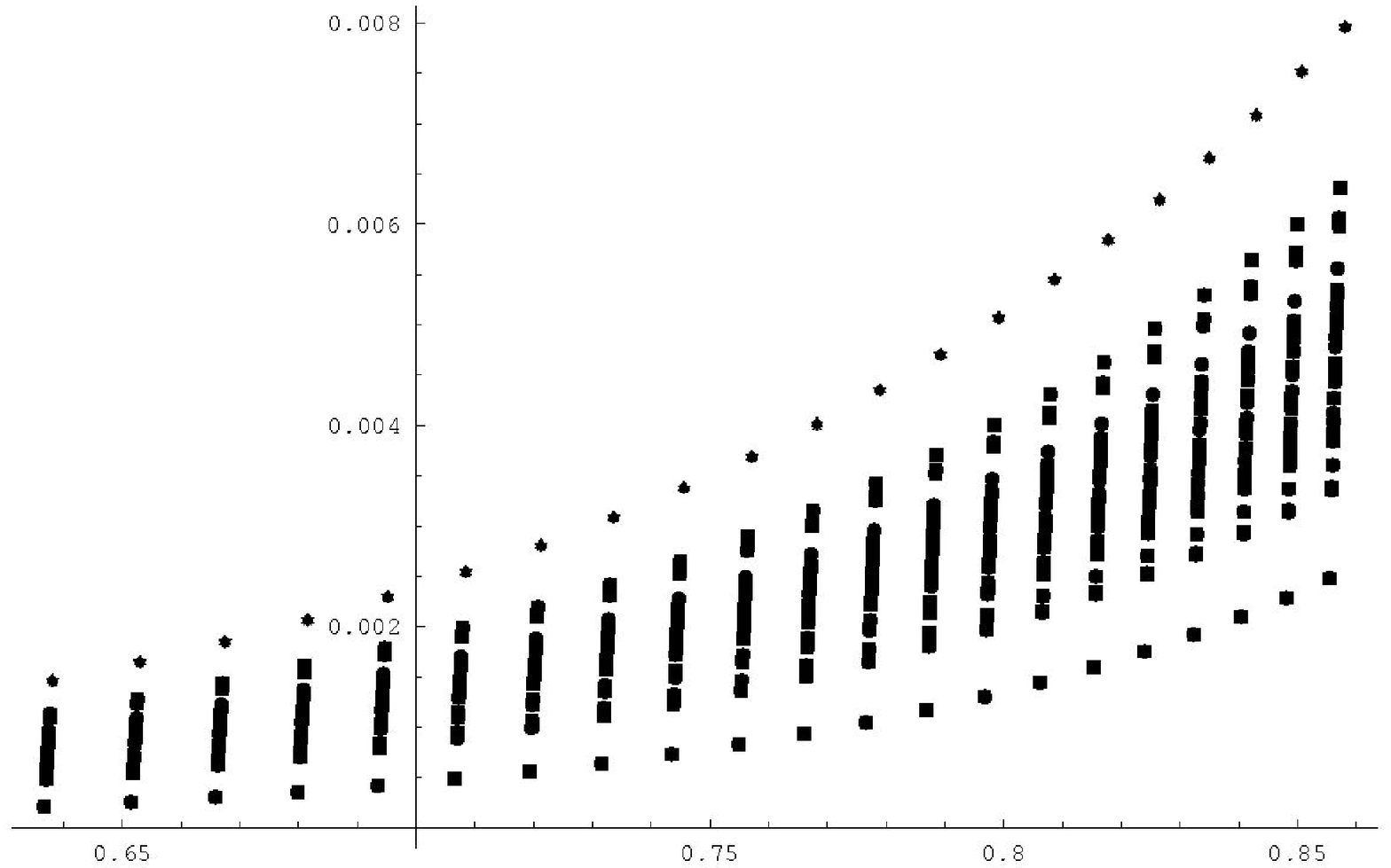}
     \caption{Difference between actual entropy and the lower
     asymptotic matching
     conjecture for 4-regular graphs with $U=C_{8}$}
     \label{fig:r4-8}
      \end{figure}

 \begin{example}[r=4]\label{exr=4}
      In our first family we let $l$, the length of the cycle
      $U=C_l$,
      vary from 4 to 8.  We tested all
      permutation matrices $A$, which give every vertex $(u,k)$ in $G_{n}$  one
      neighbor in the level $k-1$ and one in the level $k+1$, and give
      rise to a bipartite
      $G_n$. We thus have a family of bipartite 4-regular graphs which
      includes the standard square lattice tori.

      In Figure \ref{fig:r4-8} we plot the difference between
      the actual values of the entropies $h_G(p)$, for all choices of
      $A$,
      and the lower
      asymptotic matching
      conjecture for a
      given range of densities $p$. The highest curve correspond to the
      normal torus graph, i.e it is the graph with the maximum number of
      matchings of a  given size in this family.

 \end{example}

\begin{figure}[tbp]
     \includegraphics*[width=14cm,height=8cm]{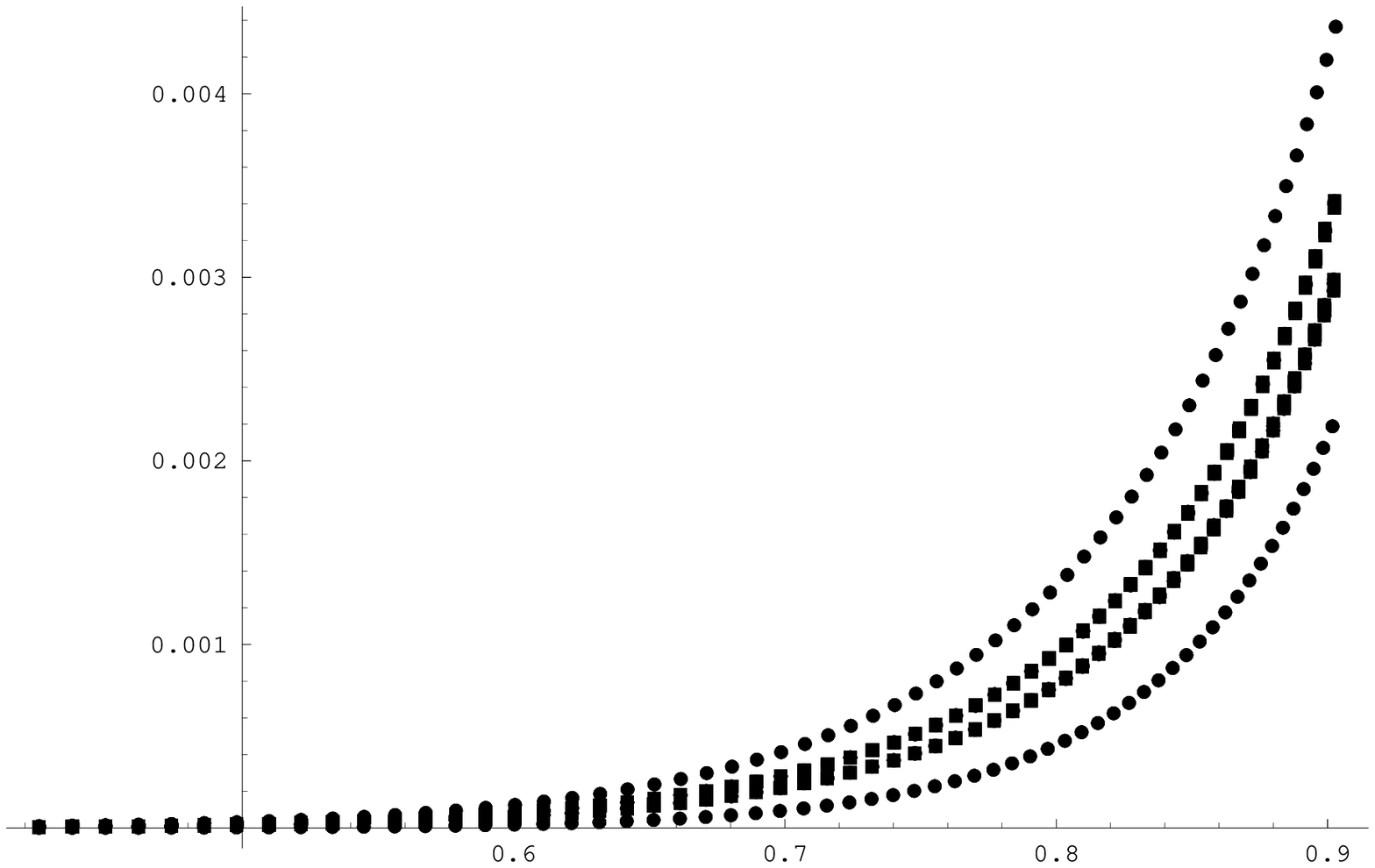}
     \caption{Difference between actual entropy and the lower
     asymptotic matching conjecture
     for 3-regular graphs with $U=C_{10}$}
     \label{fig:r3-10}
      \end{figure}

 \begin{example}[r=3]\label{exr=3}
      For our second example we again chose $U=C_l$ to be a cycle
      of length $l$. Here we
      chose $A$ so that if we number the vertices on the cycle
      $1,\ldots,l$ the even vertices in the level $k$ have an
      an odd vertex as a neighbor in the level $k+1$ for
      $k=1,\ldots,n$.  $G_n$ is cubic, and for $l$ and $n$ even
      $G_n$ is bipartite. This family includes the torus
      graphs obtained from the hexagonal lattice which all have girth at
      least 6. In this case we let the length of the cycle vary from $4$
      to 10 and again the conjecture was found to hold.
      In Figure \ref{fig:r3-10} we display
      the difference between
      the actual values of the entropies $h_G(p)$, for all choices of
      $A$, and the lower asymptotic matching
      conjecture for a given range of densities $p$.
      Here the values typically stayed closer to the conjecture
      than for the $4$-regular case, which is to be expected since this
      graph family has higher girth.

 \end{example}

 Apart from the above tests, we also tested some more arbitrarily
 chosen connection matrices giving graphs of degree 6. This was
 done by $U$ as a cycle and choosing the connection matrix $A$,
 having two $1$'s in each row and column.  Again the conjectures
 were
 found to hold but here the deviation up from the conjectured lower
 bound was even smaller. This is again expected since the
 conjecture should become more accurate for graphs of higher
 degree.


 \section{Infinite graphs with the maximal pressure}

 In this section we give a partial justification for the
 computational result in Example \ref{exr=4} that the highest curve
 correspond to the normal torus graph.

  \begin{theo}\label{maxpresex}  Let $F=(U,D)$ be an undirected
  graph and consider an infinite graph $G=(V,E)$ defined
  as follows.  $V=U\times \Z$, i.e. the vertices of $G$
  consists layers $(U,k), k\in\Z$.  The edges $E$ of $G$ connect
  two vertices on the level $k,j\in\Z$ only if $|k-j|\le 1$.
  The edges between any two vertices on level $k$ are given by
  $D$.  The edges between the level $k$ and the level $k+1$
  are given by $\#U \times \#U$ permutation matrix
  $A_k=(a_{uv,k})_{u,v\in U}$ for each $k\in\Z$.
  Thus $((u,k),(v,k+1))\in E \iff a_{uv,k}=1$.
  Let $M(t), t\in\R$ and $\tilde A_k$ be the $2^{\#U}\times 2^{\#U}$
  nonnegative matrices defined as in \S3.  Then the pressure
  $P_G(t)$ of $G$ is given as
  \begin{equation}\label{presforG}
  P_G(t)=\limsup_{i,j\to\infty}\frac{\log \tr M(t)\tilde
  A_{-j}M(t)\tilde A_{-j+1}\ldots M(t)\tilde A_i}{(i+j+1)\#U}.
  \end{equation}
  Let $G_0$ be the infinite graph obtained by letting $A_k$ to
  be the identity matrix for each $k\in\Z$.
  Then
  \begin{equation}\label{presforG0}
  P_{G_0}(t)=\frac{\log \rho(M(t))}{\#U} \ge P_{G}(t)
  \end{equation}
  for any $t\in\R$ and any $G$ of the above form.
  In particular the monomer-dimer entropy $h_G$ of $G$,
  which is equal to $P_G(0)$, does not exceed $h_{G_0}=P_{G_0}(0)$.

  \end{theo}
  \proof
  Fix $i,j\in\N$ and let $n=i+j+1$.  Define $G_n=(V_n,E_n)$ to be
  the following graph.  $V_n$ consists of $n$ copies of $U$,
  labelled as $(U,k)$ for $k=-j,\ldots,i$.  The edges of $E_n$
  are induced by the edges of $G$, except that the edges from the
  level $i$ are connected to the edges of the level $-j$, which is
  identified with the level $i+1$, by the connection matrix $A_i$.
  The arguments of \S2  yield that $\psi(e^{2t},G_n)=
  \tr M(t)\tilde
  A_{-j}M(t)\tilde A_{-j+1}\ldots M(t)\tilde A_i$.
  Hence $P_G$ is given by (\ref{presforG}) \cite{FP2}.

  Consider now the case of $G_0$ where $A_k=I$.  Then
  (\ref{presforG}) yields $ P_{G_0}(t)=\frac{\log \rho(M(t))}{\#U}$.
  (See for example \cite[\S10]{Fr3} for the self-contained details
  of the arguments on matrices used here.)  From the definition
  of $M(t)$ it follows that $M(t)$ is a nonnegative and symmetric
  matrix.  Hence $\rho(M(t))=||M(t)||$, where $||M(t)||$ is the
  $l_2$ operator norm of $M(t)$.
  Since $A_k$ is a permutation it follows that $\tilde A_k$
  is also a permutation matrix.  Hence
  $||\tilde A_k||=1$.  Thus
  \begin{eqnarray*}
  &&||M(t)\tilde A_{-j}M(t)\tilde A_{-j+1}\ldots M(t)\tilde A_i||\le
  ||M(t)||^{i+j+1}=\rho(M(t))^{i+j+1} \Rightarrow\\
  &&\tr M(t)\tilde A_{-j}M(t)\tilde A_{-j+1}\ldots M(t)\tilde
  A_i\le 2^{\#U}\rho(M(t))^{i+j+1} \Rightarrow\\
  &&P_G(t)\le \frac{\log \rho(M(t))}{\#U}.
  \end{eqnarray*}
  This proves (\ref{presforG0}).

  From the definition of monomer-dimer entropy of $G$ \cite{FP}
  it follows that $h_G=P_G(0)$.  Hence $h_G\le h_{G_0}$.
  \qed

  \begin{con}\label{maxentpex}  Let the assumptions of Theorem
  \ref{maxpresex} hold.  Then for any $p\in [0,1]$ $h_G(p)\le
  h_{G_0}(p)$.
  \end{con}


 \bibliographystyle{plain}

 \newpage


 \end{document}